\documentclass[conference]{IEEEtran}
\IEEEoverridecommandlockouts
\usepackage{cite}
\usepackage{amsmath,amssymb,amsfonts}
\usepackage{algorithmic}
\usepackage{graphicx}
\usepackage{textcomp}
\usepackage{xcolor}
\usepackage{booktabs}
\usepackage{float}
\newcommand{\ra}[1]{\renewcommand{\arraystretch}{#1}}
\def\BibTeX{{\rm B\kern-.05em{\sc i\kern-.025em b}\kern-.08em
    T\kern-.1667em\lower.7ex\hbox{E}\kern-.125emX}}
\ifCLASSOPTIONcompsoc
    \usepackage[caption=false, font=normalsize, labelfont=sf, textfont=sf]{subfig}
\else
\usepackage[caption=false, font=footnotesize]{subfig}
\fi
\usepackage{afterpage}
\begin{document}

\title{Market Value of Differentially-Private Smart Meter Data\\

\thanks{This work was supported by ESRC through  LISS DTP (ES/P000703/1) and by EPSRC through  Supergen Energy Networks (EP/S00078X/1)}
}

\author{\IEEEauthorblockN{Saurab Chhachhi}
\IEEEauthorblockA{
\textit{Imperial College London}\\
London, UK \\
saurab.chhachhi11@imperial.ac.uk}
\and
\IEEEauthorblockN{Dr. Fei Teng}
\IEEEauthorblockA{
\textit{Imperial College London}\\
London, UK \\
f.teng@imperial.ac.uk}
}

\maketitle

\begin{abstract}
This paper proposes a framework to investigate the value of sharing privacy-protected smart meter data between domestic consumers and load serving entities. The framework consists of a discounted differential privacy model to  ensure individuals cannot be identified from aggregated data, a ANN-based short-term load forecasting to quantify the impact of data availability and privacy protection on the forecasting error and an optimal procurement problem in day-ahead and balancing markets to assess the market value of the privacy-utility trade-off. The framework demonstrates that when the load profile of a consumer group differs from the system average, which is quantified using the Kullback-Leibler divergence, there is significant value in sharing smart meter data while retaining individual consumer privacy.  
\end{abstract}

\begin{IEEEkeywords}
data markets, differential privacy, load forecasting, smart grid, smart meters
\end{IEEEkeywords}

\section{Introduction}
Smart metering for domestic consumers is seen as a key enabler in moving towards a more dynamic, cost-effective, cost-reflective and decarbonised electricity network. It provides benefits for both load serving entities (LSE) and consumers through improved billing accuracy, real-time feedback on consumption and enabling innovative business models which harness demand response and dynamic pricing schemes. Access to granular data, such as half-hourly (HH) consumption of individuals, can be used to discern personal information raising issues around privacy and data usage. An overview of potential privacy concerns and smart meter data misuse is provided in \cite{Veliz2018}. Survey studies have shown that there are a range of attitudes towards smart meter data privacy. Some consumers are happy to share their consumption data, others are willing to share if details on how such data will be used and importantly how it may benefit the system as well as benefit them personally is provided, and those who are reluctant to share data under any circumstances\cite{Dickman2017}. 

Currently meter data are collected and processed on an ad-hoc basis by LSEs and then sent to the settlement body. In the UK, the introduction of smart meters will mean that, for settlement purposes, data will be collected by a centralised data company (DCC), with LSEs no longer being part of the process\cite{Ofgem2020}. Moving towards HH settlement will mean that LSEs will have to forecast HH consumption as opposed to daily volumes. This raises the question as to whether LSEs should have access to individual consumers HH data and how that access should be provided for forecasting purposes. Ofgem, the UK energy regulator, has discussed the use of privacy-preserving mechanisms but there is a lack of understanding as to the costs and benefits of such measures\cite{Ofgem2020}. 

Extant literature on smart meter privacy-preserving mechanisms has focused on their effect on data utility (i.e. change in forecasting accuracy)\cite{Hasan2020}. However the resulting impact on energy procurement costs has, to the best of our knowledge, not been investigated. To this end, we propose a framework to assess the smart meter data privacy cost-benefit trade-off for forecasting purposes, simultaneously addressing the specific privacy concerns discussed above. The framework consists of three alternative settlement and forecasting schemes: one in which HH data sharing is mandatory, one in which HH data are not shared and one where HH is shared but is privacy-preserved using differential privacy (DP). To compare the different schemes a forecasting and procurement strategy for a LSE is developed. It consists of an adaptable short-term load forecasting mechanism and an optimal procurement strategy for the LSE in the day-ahead and balancing markets. This paper makes the following contributions:
\begin{itemize}
    \item Proposes a framework to explicitly link smart meter data sharing to monetary value incorporating privacy concerns.
    \item Develops a forecasting and procurement strategy for a price-making LSE engaged in the day-ahead and balancing markets within which the privacy-utility trade-off can be assessed.
    \item Applies the Kullback-Leibler divergence as a potential indicator of data value within this context.
    \item Presents a case study using actual smart meter data and historical market prices.
\end{itemize}
In Section \ref{sec:privacy}, we review and outline privacy-preserving mechanisms for smart meter data. Section \ref{sec:framework} outlines the different load settlement schemes. The forecasting and procurement model is detailed in Section \ref{sec:setup}. Section \ref{sec:study} presents the results of a numerical case study based on actual domestic smart meter data and market prices. Finally, conclusions are drawn and future research directions discussed in Section \ref{sec:conc}.

\section{Smart Meter Data Privacy}\label{sec:privacy}
An LSE is required to forecast its consumer group's aggregated load and then purchase sufficient energy. To produce load forecasts the LSE requires historical consumption data of its consumers, which may vary in levels of aggregation or temporal resolution. As historical consumption data can provide significant amounts of personal information about an individual, the development of privacy-preserving mechanisms for releasing smart meter data has been a growing area of research. Privacy-preserving mechanisms can be categorised into the following \cite{Giaconi2018}:
\begin{itemize}
    \item Cryptographic methods such as encryption.
    \item Data manipulation which includes spatial aggregation and sampling, anonymization and differential privacy (DP).
    \item User demand shaping using batteries.
\end{itemize}
Smart meter data can be used to discern a wide range of attributes of an individual and it is difficult to specify the range and depth of potential data misuse. When defining privacy, for smart meter data, we must think about what specific information an individual wants to keep private. Although encryption, data aggregation and anonymization provide some increased privacy protection they do not guarantee that an individual cannot be identified or prevent user information leakage\cite{Hasan2020,Dwork2014}. DP offers a mechanism to ensure that an individual ($n$) cannot be specifically identified from within a dataset (e.g. an LSE's consumer group) while potentially preserving a high degree of data utility. This approach can protect an individual's privacy even as data analytics and machine learning techniques evolve and new use cases emerge, since the data cannot be attributed to the individual. DP introduces a mathematical framework to define the likelihood of being identified when making a query (in this case aggregated load) on a dataset. A data reporting mechanism is $\epsilon$-differentially private for $\epsilon>0$ if for any pair of neighbouring datasets $E(t)$,$E^{\prime}(t)$  (where the two datasets differ by only one individual) and some aggregated output the following holds\cite{Dwork2014}:
\begin{IEEEeqnarray}{C}
    \frac{p(\hat E(t)|E(t))}{p(\hat E(t)|E^{\prime}(t))} \leq exp(\epsilon)
    \label{eq:dp_def}
\end{IEEEeqnarray}
To achieve this data are obfuscated by introducing Laplacian noise to each individual load profile with 0 mean and $b$ scaling which is given by:
\begin{IEEEeqnarray}{C}
    b = \frac{\Delta f_{t}}{\epsilon}
\end{IEEEeqnarray}
where $\Delta f_{t} = \frac{\max \left(E_{n,t}\right) - \min\left(E_{n,t}\right)}{N}$, is the global sensitivity of the output (range of the individual loads in a given period $t$), $\epsilon$ is the privacy budget which indicates the risk of being identified and $N$ is the total number of individuals in the dataset.

Most extant literature on the application of DP and its variants to smart meter data have assumed that the privacy budget is fixed and that queries are independent. A detailed review can be found in \cite{Hasan2020}. However, smart meter data and the resulting queries are continuously generated and updated. As a result the privacy loss defined by DP is accumulated across each query, requiring the addition of increasing amount of noise to ensure \eqref{eq:dp_def} holds. Over time this degrades data quality and renders new data useless\cite{Farokhi2020}. Techniques have been proposed to overcome this limitation such as selective sampling based on time series dynamics\cite{Kim2018} which improve performance but are still sensitive to the number of queries and would not guarantee the specified privacy budget over an infinite time horizon. \cite{Farokhi2020} overcomes this, providing a bounded mechanism by introducing the notion of discounted differential privacy (DDP). It draws upon the concept of discounting from economic theory to propose that data further from the past is less sensitive than current data. The resulting noise scaling can be modelled as a function of the privacy budget (as before, $\epsilon$) and the discount rate (a measure of how much one values past data, $\gamma$):
\begin{IEEEeqnarray}{C}
    b = \frac{\Delta f_{t}}{\epsilon \left(1-\gamma\right)}, \gamma \in [0,1)
\end{IEEEeqnarray}
If one does not place any value on the privacy of past data then $\gamma = 0$ and the mechanism is equivalent to DP whereas if one places high value on the privacy of past data then $\gamma \to 1$ and the required noise tends to infinity.

\section{Domestic Load Forecasting and Settlement}\label{sec:framework}
\subsection{Non-Half-Hourly Settlement (NHHS)}
In the absence of HH smart meter data, electricity settlement is based on system-wide daily load coefficients ($DLC$) which are published ahead of time. DLCs are standardised load profiles which specify the amount of annual consumption a specific consumer group (domestic, SME etc.) consumes in a particular half hour. These are generated based on HH measurement taken from a sample of consumers within each consumer group (for details see \cite{Elexon2018}). An LSE is only required to forecast daily demand ($E^{d}$) and is therefore insulated from HH changes while still exposed to HH prices (see Table \ref{tab:mechanism}).
\begin{table}[htbp]
    \caption{Settlement Schemes}
    \begin{center}
        \resizebox{0.45\textwidth}{!}{
        \ra{1.5}
        \begin{tabular}{c c c c c}
        \toprule
        & \textbf{NHHS}&\textbf{HHS - $DLC^{sys}$}& \textbf{HHS - $E^{hh}$}& \textbf{HHS – $DDP(\epsilon,\alpha)$}\\
        \midrule
        \textbf{Forecast Input}&$E^{d}$, $DLC^{sys}$&$E^{d}$, $DLC^{sys}$&$E^{hh}$&$E^{hh} + Lap(\epsilon, \alpha)$\\
        \textbf{Forecast Parameter}&$E^{d}$&$E^{hh}$&$E^{hh}$&$E^{hh}$\\
        \textbf{Settlement}&$E^{d}DLC^{sys}$&$E^{hh}$&$E^{hh}$&$E^{hh}$\\
        \textbf{Risk Exposure}& $E^{d}, \lambda^{bal}$ &$E^{hh}, \lambda^{bal}$&$E^{hh}, \lambda^{bal}$&$E^{hh}, \lambda^{bal}$\\
        \bottomrule
        \end{tabular}}
    \label{tab:mechanism}
    \end{center}
\end{table}
\subsection{Half-Hourly Settlement (HHS)}
To assess what the underlying value of sharing HH data would be we present three alternatives (see Table \ref{tab:mechanism}): a scheme in which data sharing is mandatory i.e. the LSE has access to all its consumers aggregate unaltered HH data (HHS - $E^{hh}$), a scheme where only aggregate daily data are shared (HHS - $DLC^{sys}$) and a scheme where aggregate HH data are shared but is privacy-protected using DDP (HHS - $DDP(\epsilon, \gamma)$). Under all these schemes settlement is based on actual HH consumption but the data available for forecasting purposes differs. An overview of the dataflows for each scheme is shown in Fig. \ref{fig:framework}. The next section details the forecasting and procurement models used.

\afterpage{
\begin{figure}[htbp]
\centerline{\includegraphics[width=0.45\textwidth]{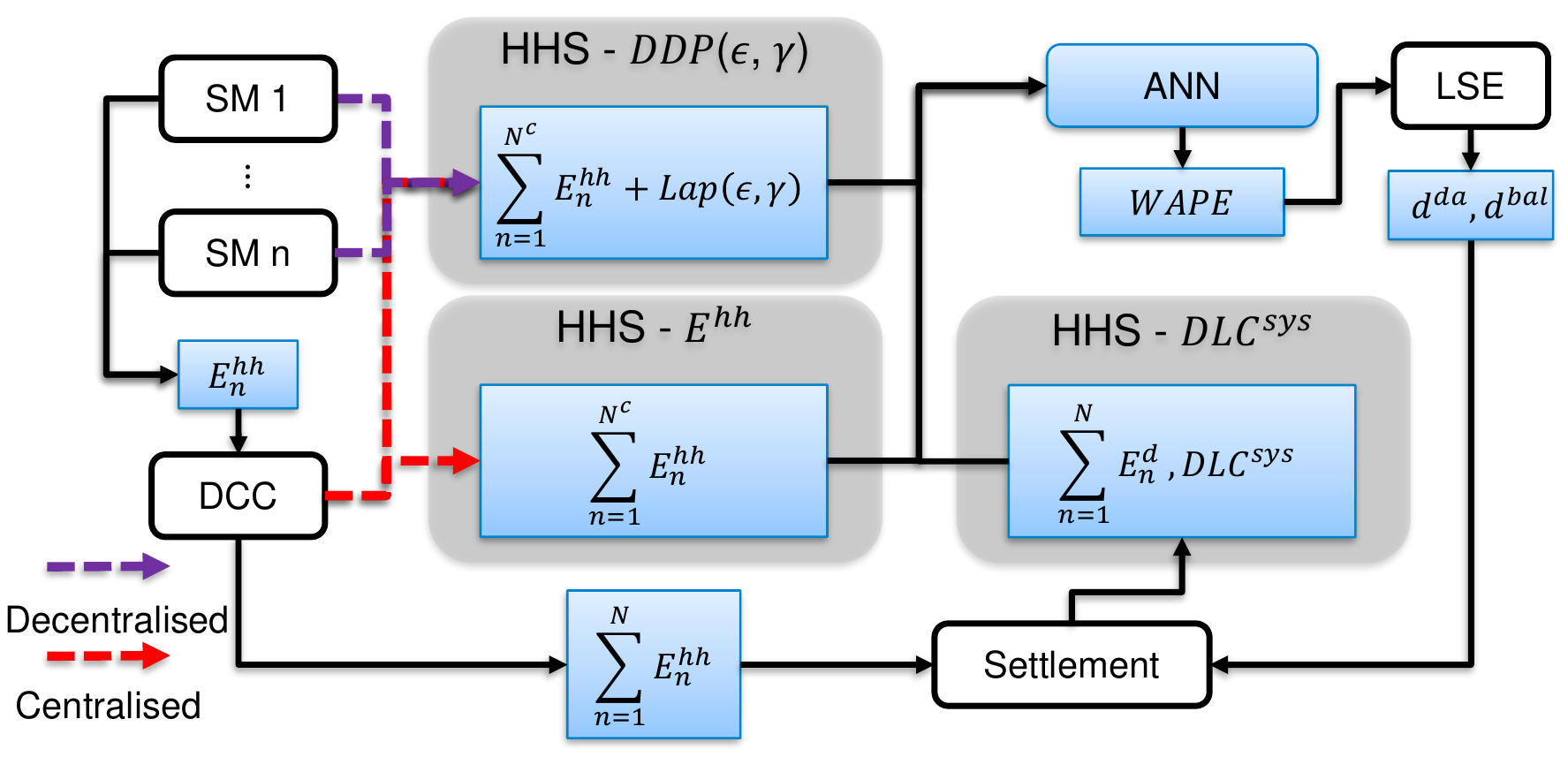}}
\caption[Caption for LOF]{Overview of Proposed Framework.\footnotemark}
\label{fig:framework}
\end{figure}
\footnotetext{Laplacian noise can be constructed by summing $n$ i.i.d. Gamma distributions allowing for decentralised noise addition at the smart meter\cite{Acs2011}.}
}

\section{Model Definition}\label{sec:setup}
\subsection{Short-Term Load Forecasting}\label{sec:fore}
Artificial Neural Networks (ANN) have been widely used and perform well for short-term forecasting applications and are able to capture both linear and non-linear dependencies\cite{Wang2018, Vagropoulos2015}. We use a simple ANN consisting of three layers: input layer, hidden layer, and output layer where the hidden layer has four hidden neurons. The following features are considered:
\begin{align}
    \begin{split}
        X = [W, WD, SP, E_{t-h}, E_{t-h-1}, E_{t-2h+1},\\
        E_{t-2h}, E_{t-3h}]
    \end{split}
\end{align}
where $E_{t-^{*}h}$ are the lagged/historical load values and h is the number of periods in the day (48), $W$ is the week in the year, $WD$ is the day of the week, and $SP$ is the settlement period. The model is implemented in Python using the MLPRegressor model in Scikit-Learn.
\subsection{Load Serving Entity (LSE) Procurement Problem}
A LSE needs to procure energy to meets its customer group's load by participating in long-term trading, day-ahead and intra-day markets, and settling any imbalances between purchase volumes and actual consumption in the balancing market. In this paper we focus on day-ahead and balancing markets. The LSE's procurement strategy can be formulated as a two-stage risk-constrained stochastic program similar to \cite{Song2018}. We assume the LSE is a price-making market entity in both the day-ahead and balancing market and account for risk-aversion by including the optimisation of the conditional value-at-risk (CVaR). The formulation is as follows:
\begin{subequations}\label{eq:opti}
\begin{align}
    \begin{split}
        \min_{d^{da}, d^{bal}} \overbrace{\underbrace{\sum_{t}\lambda^{da}_{t} d^{da}_{t}}_{Day-Ahead} + \underbrace{\sum_{s}\pi_{s}\sum_{t}\lambda^{bal}_{s,t} d^{bal}_{s,t}}_{Balancing}}^{\text{Expected Cost} (\hat\Omega)}\\
        + \underbrace{\beta \left[\zeta + \frac{1}{1-\alpha} \sum_{s} \pi_{s}\eta_{s} \right]}_{CVaR}
    \end{split}
    \tag{\ref{eq:opti}}
\end{align}
s.t.
\begin{IEEEeqnarray}{c}
     d^{da}_{t} + d^{bal}_{s,t} = d^{fore}_{t} + d^{err}_{s,t}, \forall s, \forall t \\
     \sum_{t}\left[\lambda^{da}_{t} d^{da}_{t} + \sum_{t}\lambda^{bal}_{s,t} d^{bal}_{s,t}\right] - \zeta \leq \eta_{s}, \forall s \\
    \eta_{s} \geq 0, \forall s
\end{IEEEeqnarray}
where $d^{da}_{t}$ and $d^{bal}_{s,t}$ are the volumes procured by the LSE in the day-ahead and balancing market respectively, $\lambda^{da}_{t}$ and $\lambda^{bal}_{s,t}$ are the market prices, $\pi_{s}$ is the scenario probability, $\beta$ is the risk-aversion factor, $\alpha$ is the CVaR confidence interval, $\zeta$ and $\eta_{s}$ are auxiliary variables to calculate CVaR, $d^{fore}_{t}$ and $d^{err}_{s,t}$ are the day-ahead forecast load and the realised error. Given that the LSE is a price-making entity $\lambda^{da}_{t}$ and $\lambda^{bal}_{s,t}$ are dependent on the total demand. These can be modelled as piece-wise linear curves:
\begin{IEEEeqnarray}{c}
    \lambda^{da}_{t} = \sum_{b} \lambda^{G}_{b} u^{da}_{t,b}, \forall t \\
    \lambda^{bal}_{s,t} = \sum_{f} \lambda^{F}_{b} u^{bal}_{s,t,f}, \forall s, \forall t \\
    D^{sys}_{t} - \frac{\Delta}{2} \leq \sum_{b} u^{da}_{t,b} \Tilde{D}^{sys}_{b} \leq D^{sys}_{t} + \frac{\Delta}{2}, \forall t \\
    D^{imb}_{s,t} - \frac{\Delta}{2} \leq \sum_{f} u^{bal}_{s,t,f} \Tilde{D}^{imb}_{f} \leq D^{imb}_{s,t,} + \frac{\Delta}{2}, \forall s, \forall t
\end{IEEEeqnarray}
where $D^{sys}_{t} = D^{da}_{t} + d^{da}_{t}$, the total system demand in period $t$ and $D^{imb}_{s,t} = D^{bal}_{s,t} + d^{bal}_{s,t}$, the total system imbalance in scenario $s$, $\Tilde{D}^{da}_{b}$ is a discretisation of the system demand into $b$ increments of $\Delta$ (similarly for $\Tilde{D}^{bal}_{f}$) and $u^{da}_{t,b}$ and $u^{bal}_{s,t,f}$ are binary variables which select the appropriate demand level. The resulting model is a MIQP due to the products of binary and continuous variables ($u^{da}_{t,b} d^{da}_{t}$ and $u^{bal}_{s,t,f} d^{bal}_{s,t}$). These can be linearised by replacing the bilinear terms with a new variable on which a number of constraints are imposed giving an exact MILP reformulation\cite{williams2013}. For example the term $u^{da}_{t,b} d^{da}_{t}$ can replaced by an auxiliary variable, $C^{da}_{t,b}$, and four additional constraints:
\begin{IEEEeqnarray}{c}
    C^{da}_{t,b} \leq u^{da}_{t,b} d^{max}_{t,b}, \forall t, \forall b \\
    C^{da}_{t,b} \leq d^{da}_{t}, \forall t, \forall b \\
    C^{da}_{t,b} \geq d^{da}_{t} - (1-u^{da}_{t,b}) d^{max}_{t,b}, \forall t, \forall b \\
    C^{da}_{t,b} \geq 0, \forall t, \forall b 
\end{IEEEeqnarray}
\end{subequations}
The MILP reformulation is implemented in FICO\textsuperscript{TM} Xpress through the Python API.
\subsection{Assessment Metrics}
To gauge the difference between the load profile of the LSE's consumer group and the rest of the system we employ the Kullback-Leibler divergence (KLD). Its value can be interpreted as the information gain achieved if HH data from the consumer group ($DLC^{c}$) is used instead of system level data ($DLC^{sys}$). We assume the average weekly $DLC$ variations are normally distributed. In this context it can be defined as follows\cite{Roberts2002}:
\begin{IEEEeqnarray}{C}
    KLD = \sum_{t \in T^{w}} \left[\log \frac{\sigma^{sys}_{t}}{\sigma^{c}_{t}} +
         \frac{\sigma^{c^{2}}_{t} + \left(\mu^{sys}_{t} - \mu^{c}_{t}\right)^{2}}{2\sigma^{sys^{2}}_{t}} - \frac{1}{2} \right]
\end{IEEEeqnarray}
where $\mu$ and $\sigma$ are the mean and standard deviation of the $DLC$ for each HH of the week respectively.

To measure the accuracy of the forecasts we employ the weighted-absolute percentage error (WAPE) metric as it exhibits more stable behaviour for values close to zero. This is especially relevant in the presence of PV and battery storage as net load profiles can be negative.
\begin{IEEEeqnarray}{C}
    WAPE = \frac{\sum_{t} \left\lvert E^{act}_{t} - E^{fore}_{t} \right\rvert}{\sum_{t} E^{act}_{t}}
\end{IEEEeqnarray}

\section{Case Study}\label{sec:study}
\subsection{Data}
We use smart meter data from the CER Behavioural Trials which includes 6010 residential and SME consumers for a period of 75 weeks\cite{CER2012}\footnote{After data processing to remove periods and meters for which less than 95\% of data was available.}. For 50\% of consumers synthetic PV \cite{Pfenninger2016} and EV \cite{element} load profiles are added to better reflect the increased load diversity expected in the future. Day-ahead and balancing market bidding curves are generated based on historical UK market data for 2018 from Elexon. The WAPEs of the ANN described in Section \ref{sec:fore} are used to 50 generate demand forecast scenarios, assuming they are normally distributed and then scaled to a representative UK system level based on the share of total load of consumers in the CER dataset.

\begin{figure}[htbp]
    \centering
  \subfloat[KLD for Different Market Shares\label{fig:baseline_a}]{
       \includegraphics[width=0.47\linewidth]{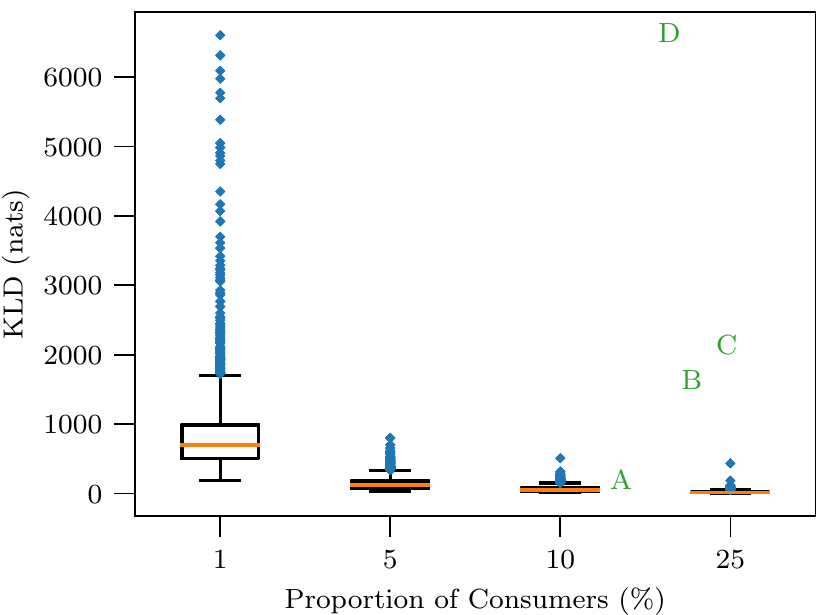}}
    \hfil
  \subfloat[WAPE for Select KLD\label{fig:baseline_b}]{
        \includegraphics[width=0.47\linewidth]{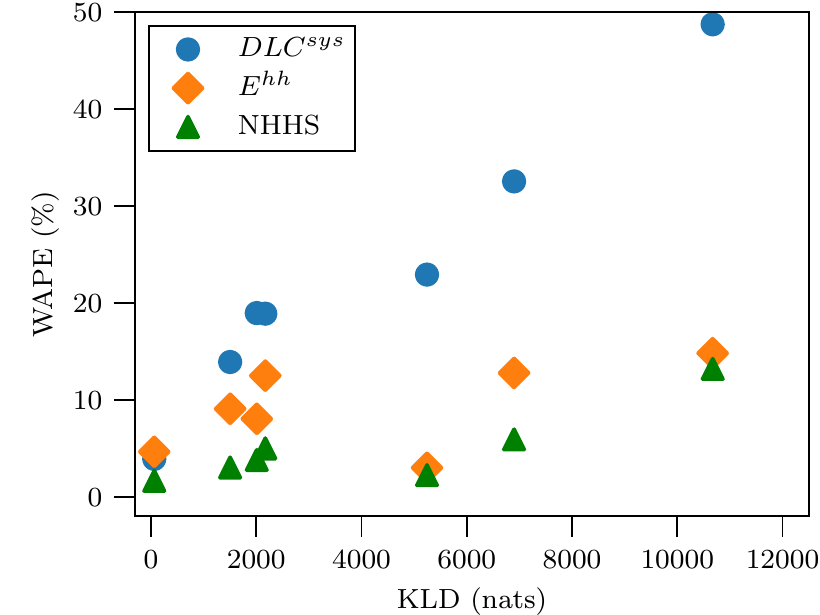}}
  \caption{KLD and WAPE for 2000 Sampled Groups.}
  \label{fig:baseline} 
\end{figure}

Fig. \ref{fig:baseline_a} shows the KLDs for randomly sampled meters under various proportions of consumers. When the meters are a small proportion of the total consumers the KLDs of the aggregated load can be large but as the proportion increases the KLDs decrease significantly. We argue that as LSEs begin to offer more innovative and targeted tariff mechanisms KLDs could be large even for large groups of consumers as they change consumption patterns based on particular time-varying incentives. We select four groups based on K-Means clustering of average DLC for each consumer across the week (shown in grey in Fig. \ref{fig:clusters}) to test the framework. To add context we plot the resulting clusters on Fig. \ref{fig:baseline_a} (A-D). From Fig. \ref{fig:baseline_b}, which shows the forecast error, it is clear that as the KLD of a consumer group increases, having HH data for that specific group results in a greater reduction in forecasting error.

\subsection{Results}
\subsubsection{Forecast Accuracy for Different Consumer Group}
\begin{figure*}[htbp]
    \centering
    \includegraphics[width=0.24\linewidth]{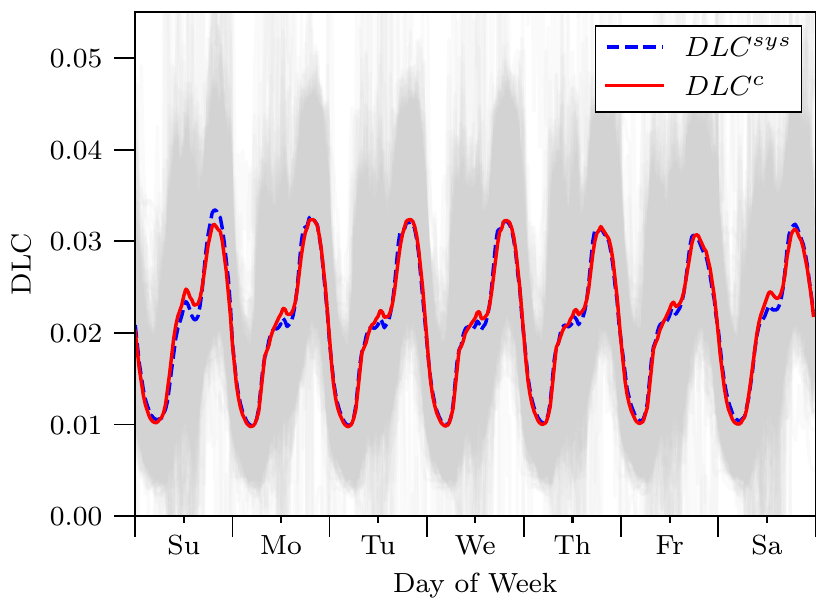}
    \hfil
    \includegraphics[width=0.24\linewidth]{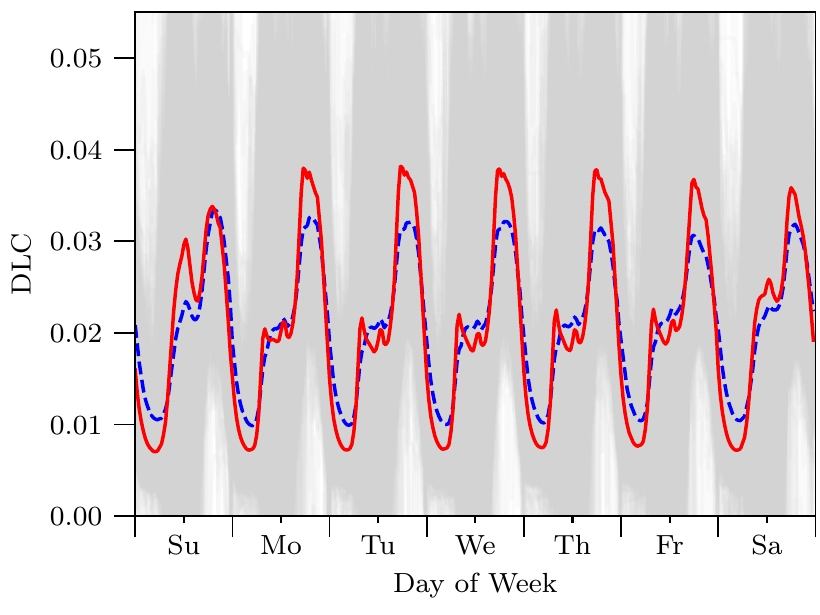}
    \hfil
    \includegraphics[width=0.24\linewidth]{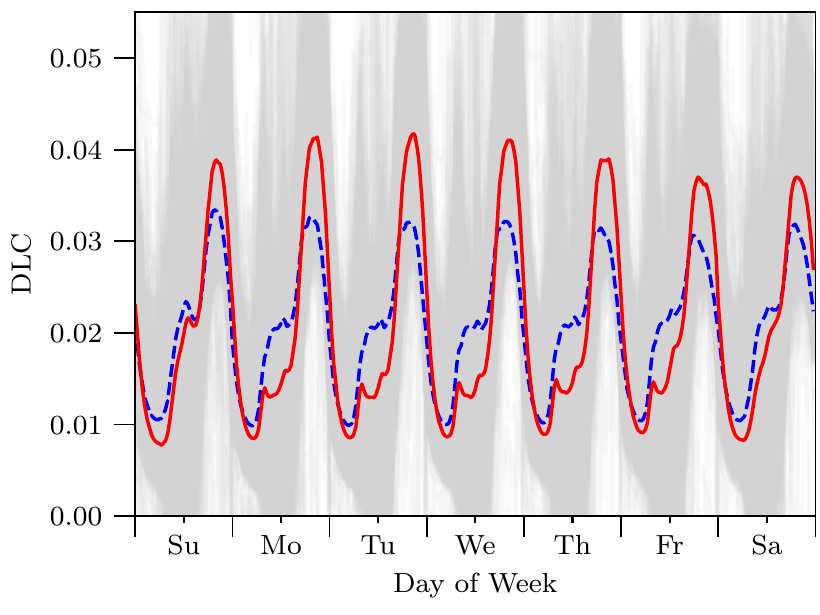}
    \hfil
    \includegraphics[width=0.24\linewidth]{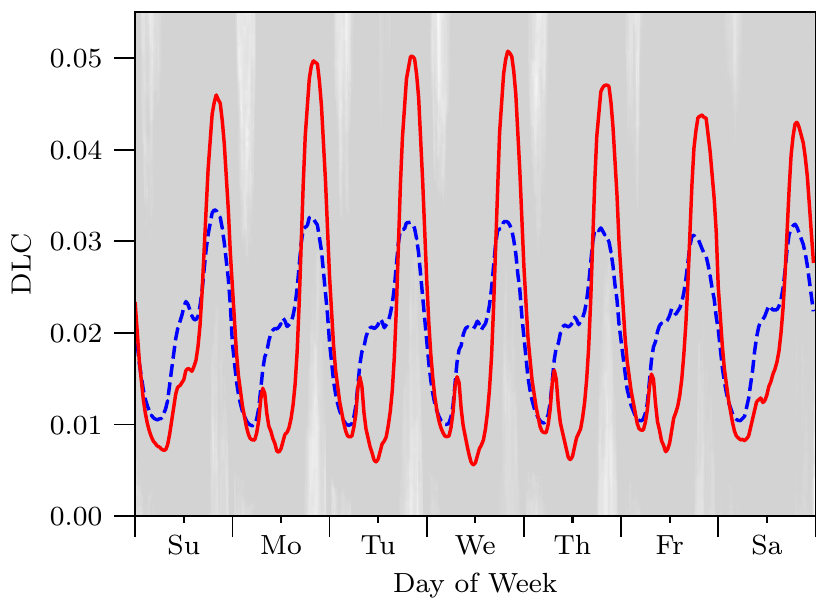}\\
    \subfloat[Group A ($KLD = 62$)\label{fig:effect_DDP0}]{
        \includegraphics[width=0.24\linewidth]{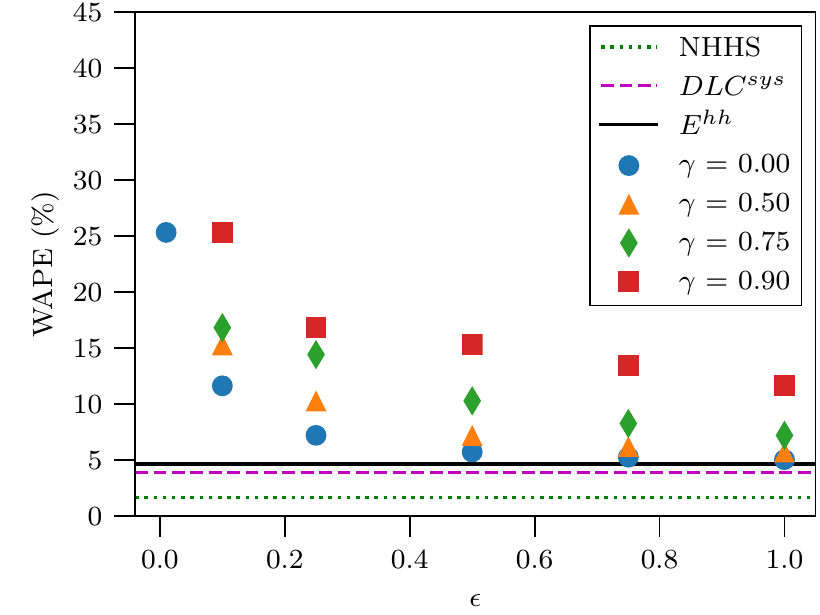}}
        \hfil
    \subfloat[Group B ($KLD = 1504$)\label{fig:effect_DDP1}]{
        \includegraphics[width=0.24\linewidth]{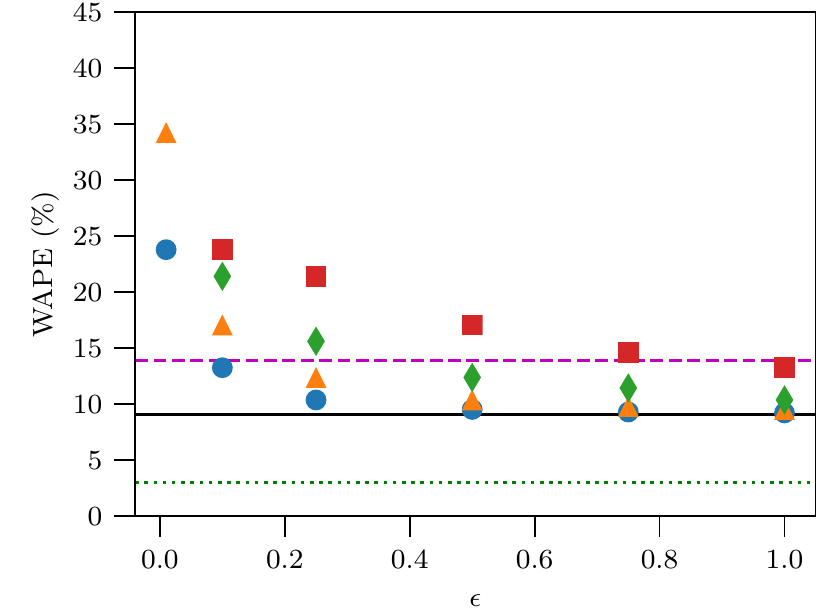}}
        \hfil
    \subfloat[Group C ($KLD = 2010$)\label{fig:effect_DDP2}]{
        \includegraphics[width=0.24\linewidth]{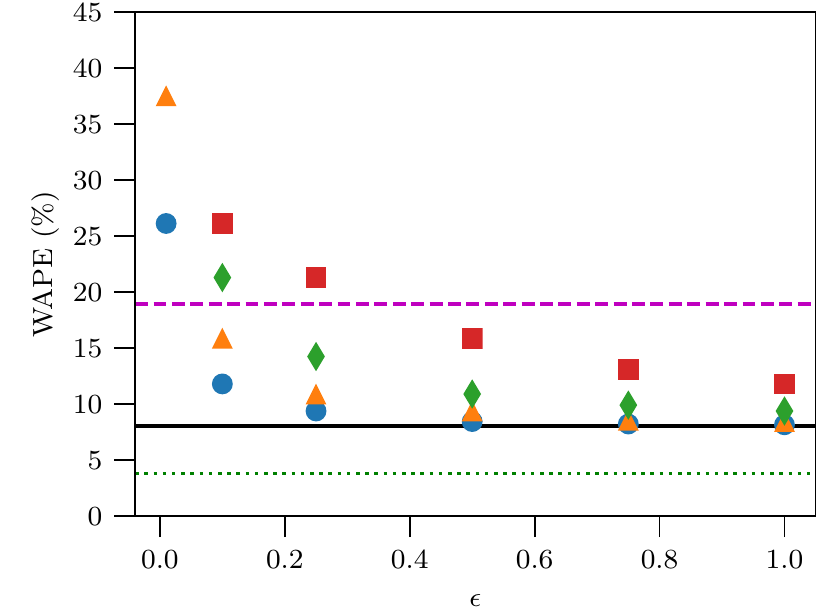}}
        \hfil
    \subfloat[Group D ($KLD = 6899$)\label{fig:effect_DDP3}]{
        \includegraphics[width=0.24\linewidth]{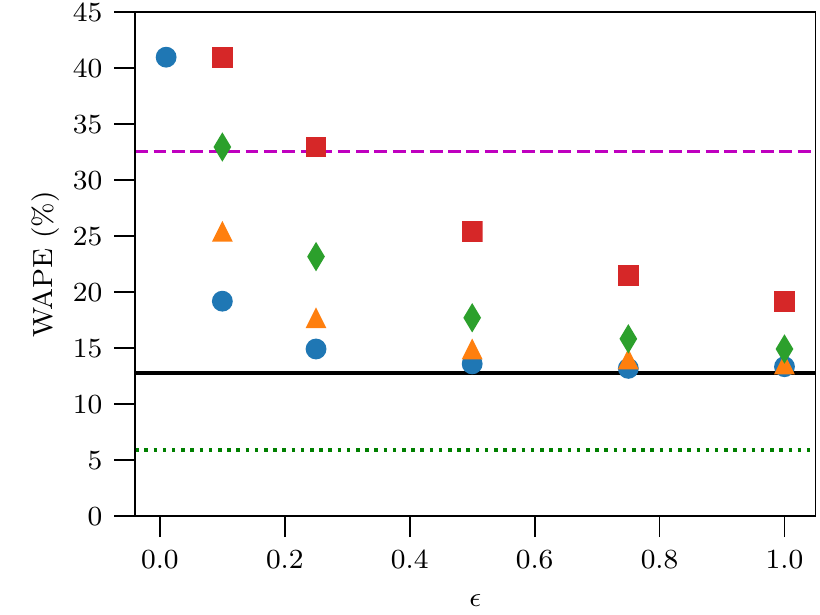}}
    \caption{Group Weekly Average DLC (top), Forecasting Error for Different Settlement Mechanisms (bottom).}
    \label{fig:clusters}
\end{figure*}

The top plots in Fig. \ref{fig:clusters} show the weekly average DLCs for the selected consumer groups. The bottom plots show the WAPE under each scheme. It is clear that when the average DLC of the consumer group is similar to the system (Group A), reflected in a low KLD, the WAPE is low even in the case where HH data are not provided ($DLC^{sys}$). As a result providing HH data using the DDP mechanism does not increase utility. However, as the group KLD increases, there is an increase in the difference between the WAPE with HH data and the WAPE without access to HH data. This provides a range within which an LSE is able to explore the privacy-utility trade-off. For example for Group C the LSE would be able to gain a 5\% reduction in WAPE while providing a $(\epsilon = .25,\gamma = .75)$ level of privacy. This shows that privacy-preservation can be achieved without significantly degrading data utility.
\subsubsection{Market Value}
Fig. \ref{fig:mval} shows the procurement costs based on scenarios generated for the different settlement schemes for Group C. An LSE can make significant cost savings while still providing consumers with privacy. On average we see that a 1\% increase in WAPE results in a 0.8 - 1\% cost reduction. A greater reduction is observed in the CVaR as a 1\% increase in WAPE results in a 2-3\% reduction in CVaR. The value of better forecasting accuracy and hence HH data are also highly dependent on the market dynamics. At peak times uncertainty is more  expensive, as there is less flexibility in the balancing market when overall demand is high, resulting in larger cost differences between the schemes.
\subsubsection{Heterogeneous Privacy Preferences}
\begin{figure}[H]
    \centering
    \subfloat[Procurement Cost \label{fig:mval}]{
        \includegraphics[width=0.47\linewidth]{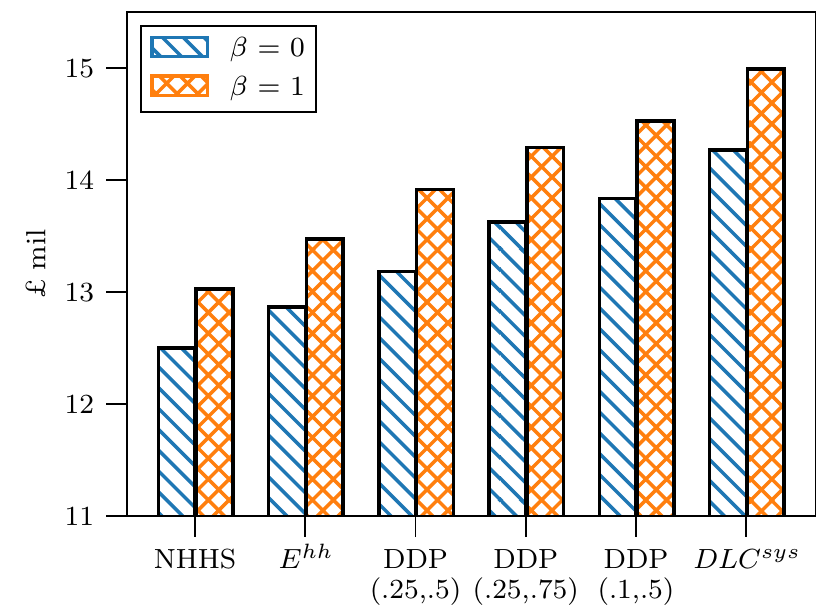}}
    \hfil
    \subfloat[Effect of Heterogeneity\label{fig:hetero}]{
        \includegraphics[width=0.47\linewidth]{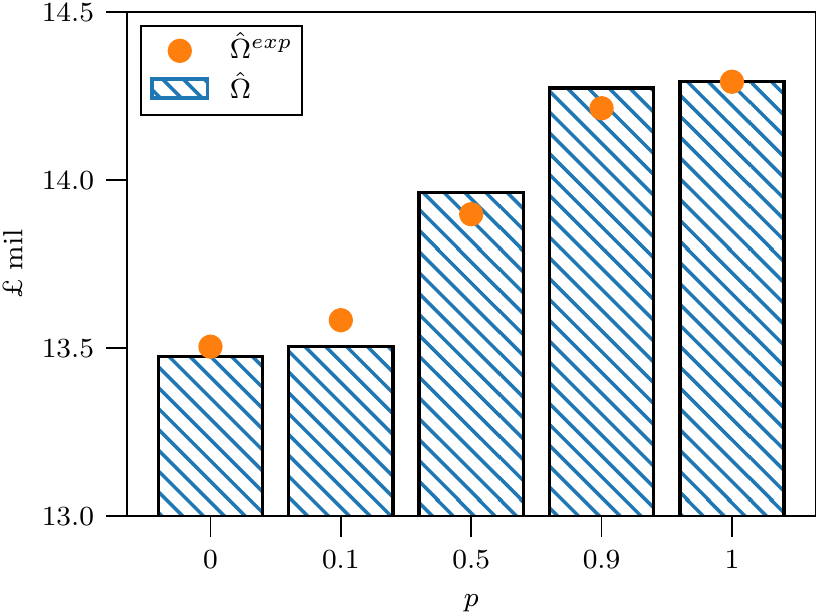}}
    \caption{Procurement Cost for the Different Schemes.}
    \label{fig:market_value}
\end{figure}

As privacy concerns vary, we investigate how the costs change when only a fraction of the consumer group has privacy concerns. Assuming a proportion of the consumer group $p$ has privacy concerns we generate forecasts separately for them using $DDP$, with the load for the remaining consumers modelled using $E^{hh}$. Fig. \ref{fig:hetero} shows the resulting procurement costs using this method for $DDP(.25,.75)$. The dots represent the weighted average cost ($\hat \Omega^{exp}$) that would be expected based on the proportion $p$. Splitting consumers based on privacy concerns can improve overall data quality and reduce overall procurement costs when $p$ is low. However a trade-off is observed between reduced data degradation, as less noise is added, and benefits of aggregation, which smoothens the load profile.
\section{Conclusions and Future Work}\label{sec:conc}
This paper investigated the value of sharing smart meter data applying a framework consisting of a discounted differential privacy model to ensure individuals cannot be identified from aggregated data, a short-term load forecasting method using ANN to quantify the impact of data availability and privacy protection on the forecasting error,  and an optimal procurement problem, to assess the market value of the privacy-utility trade-off introduced by DDP. Results show that when the load profile of a LSE's consumer group differs from the system average, which is increasingly relevant with the introduction of dynamic tariffs, and distributed storage and generation, there is significant value in sharing data while retaining individual consumer privacy. Further work is needed to assess how the benefits of smart meter data sharing can be distributed through the development of privacy differentiated tariffs or data markets to incentivise data sharing. In addition, reducing the global sensitivity parameter ($\Delta f$) by optimising the noise introduced by the DDP mechanism and explicitly incorporating heterogeneous privacy preferences would further improve performance. This framework could also be extended to forecasting demand response, net load and flexibility as well as including additional data streams such as user preferences and appliance information.
\bibliographystyle{IEEEtran}
\bibliography{IEEEabrv,reference}
\end{document}